 

\baselineskip=14pt
\parskip=10pt
\def\Tilde{\char126\relax}
\def\halmos{\hbox{\vrule height0.15cm width0.01cm\vbox{\hrule height
 0.01cm width0.2cm \vskip0.15cm \hrule height 0.01cm width0.2cm}\vrule
 height0.15cm width 0.01cm}}
\font\eightrm=cmr8  
\font\eighttt=cmtt8
\magnification=\magstephalf

\parindent=0pt
\overfullrule=0in
 
\bf
\centerline
{A Heterosexual Mehler Formula for the Straight
Hermite Polynomials (\`A La Foata)}
\rm
\bigskip
\centerline{ {\it Doron ZEILBERGER}\footnote{$^1$}
{\eightrm  \raggedright
Department of Mathematics, Temple University,
Philadelphia, PA 19122, USA. 
{\eighttt zeilberg@math.temple.edu} \hfill \break
{\eighttt http://www.math.temple.edu/\Tilde zeilberg/   .}
July 13, 1998.  
} 
}
 
{\bf Abstract:} The celebrated Foata combinatorial model for
Hermite polynomials, and his seminal and beautiful proof
of the Mehler formula, are straightened to deal with
two sexes rather than one, with the exclusion of same-sex
relationships (both marital and non-marital).
 
{\bf Introduction} 
 
In my invited talk `The unreasonable effectiveness of Combinatorics
in Orthogonal Polynomials I', 
(July 1, 1998, 10:30-11:30 local time),
at the IWOP '98, L\'eganes, Madrid,
(organized by Paco Marcellan and Renato Alvarez-Nodarse), I described,
for most of the hour, Dominique Foata's[F1] gorgeous and seminal
proof of Mehler's formula. As I got carried away in the exposition,
I replaced the original abstract `vertices' 
and `edges' by `people' and `relationships' respectively.
 
In Foata's proof
edges may have two different colors. In my more human rendition,
I called the colors of the relationships `marriage' and
`affair' (lover-mistress), which I found more engaging.
Alas, as I was adapting Foata's language to my own's,
I realized that the analogy is not complete. Traditional
Marriage and even traditional Romance assumes the existence of
two sexes, and the vertices connecting an edge to be of opposite sex.
I quickly recovered, and said that the model is liberated and modern,
and (almost) anything goes, i.e. gender is irrelevant, but
each person is allowed at most one spouse and at most one lover.
 
While I have nothing against same-sex unions, it is still
an interesting combinatorial problem to have a two-sex
analog of the Hermite polynomials, and a heterosexual
Mehler formula. It is a testimony to the power of Foata's
combinatorial approach that the adaptation is natural
and straightforward, and moreover leads to an even deeper
(at least from the WZ-point of view) Mehler-type formula.
 
{\bf The Two-Sex (Bipartite) Exponential Formula}
 
The Exponential formula ([FS][BG][F2], see [W] for a great exposition),
is one of the pillars of modern enumeration. Let's briefly 
adapt it to the bipartite case.
 
Let $A_{m,n}$ be the weight-enumerator of {\it labelled} structures
involving $m$ men and $n$ women. Let $B_{m,n}$ be the corresponding
quantity for another structure. Then if $C_{m,n}$ is the
weight-enumerator of the composite {\it labelled}
structures of {\it ordered} pairs $(a,b)$ where $a$ is an $A$-structure and 
$b$ is a  $B$-structure then, of course
$$
C_{m,n}=\sum_{k=0}^m\sum_{s=0}^n {{m} \choose {k}}{{n} \choose {s}}
A_{k,s}B_{m-k,n-s} \quad ,
\eqno(PAIRS)
$$
since among the $m$ male- and $n$
female- available  $C$- labels one has to 
choose which ones will be taken by the $A$- structure, and which
ones by the $B$- structure.
 
Now define the {\it Bivariate Exponential Generating Function}
(BiEGF) of a $2$-sex labelled structure $A$ by:
$$
f_A(t,s):=\sum_{m=0}^{\infty}\sum_{n=0}^{\infty}
{{A_{m,n}t^ms^n} \over {m!n!}} \quad.
$$
Equation $(PAIRS)$ implies
$$
f_C(t,s)=f_A(t,s)f_B(t,s) \quad .
$$
 
If $A$ and $B$ coincide, then we have that the BiEGF of
{\it ordered} pairs of labelled structures $(a_1,a_2)$, 
is $f_A(t,s)^2$. By iterating, the BiEGF of {\it ordered}
$k$-tuples is $f_A(t,s)^k$. Hence the BiEGF of {\it unordered}
$k$-tuples is $f_A(t,s)^k/k!$, and the BiEGF of all
{\it unordered} tuples is:
$$
\sum_{k=0}^{\infty} {{f_A(t,s)^k} \over {k!}}
\,=\, \exp \, [\,f_A(t,s) \, ]  \quad.
$$
 
{\bf The Straight Hermite Polynomials}
 
Let there be $m$ labelled men and $n$ labelled women.
Any person may choose to either marry or remain single.
A {\it marital profile} consists of all these choices being
made. The weight of each married couple is $x$, and hence if a
marital profile has $k$ couples, then its weight is $x^k$.
Let $H_{m,n}(x)$ be the sum of all the weights of all possible
marital profiles. Now,
$$
H_{m,n}(x)=\sum_{k=0}^{min(m,n)} 
{{m} \choose {k}}{{n} \choose {k}}k! x^k \quad,
$$
since there are ${{m} \choose {k}}$ ways of choosing which
$k$ men will marry and 
${{n} \choose {k}}$ ways of choosing which
$k$ women will marry, and there are $k!$ ways of matching them
up. 
 
The ``connected components'' of this structure are a single
man (BiEGF=$t$), a single woman (BiEGF=$s$), and a married
couple (BiEGF=$xts$), hence the BiEGF for
the Straight Hermite polynomials is:
$$
\sum_{m=0}^{\infty}\sum_{n=0}^{\infty}
{{H_{m,n}(x)t^ms^n} \over {m!n!}}=\, e^{t+s+xts} \quad.
$$
 
{\bf Legalizing Affairs}
 
Now assume that it is legal, in addition to having the option to have one
spouse (of the opposite sex), to also have the option of having one lover
(of the opposite sex). So each person may opt to have
no mates at all, to only have a spouse, to only have
a lover, or to have both. 
Multiple spouses and/or multiple lovers are forbidden.
Let's call a {\it marital-extramarital profile} such a choice.
And let the weight of an affair be $y$. Then if there are
$k$ married couples and $l$ pairs of lovers, then the
weight of the profile is $x^k y^l$. Of course, the
sum of the weights of all such profiles, when $m$ men
and $n$ women are present is
$H_{m,n}(x)H_{m,n}(y)$, since the decision whether and whom to marry
is independent of the decision whether and whom to love. Note that
it is perfectly legal to have your spouse and lover coincide.
 
{\bf The Heterosexual Mehler Formula}
 
The BiEGF of marital-extramarital profiles is
$$
\sum_{m=0}^{\infty}\sum_{n=0}^{\infty}
{{H_{m,n}(x) H_{m,n}(y) t^ms^n} \over {m!n!}}\quad .
$$
 
Let's use the Two-Sex Exponential Formula
to derive an alternative expression for the above BiEGF, by analyzing
the possible {\it connected components} of the marital-extramarital
structure.
 
Case I: A single unmarried and unloved man: BiEGF=$t$.
 
Case Ia: A single unmarried and unloved woman: BiEGF=$s$.
 
Case II: Boris is married to Natasha, Natasha loves Sasha,
Sasha is married to Lena, Lena loves Misha, Misha is married to
Sofia, ...,
 
a chain that starts with a married and unloved
man and ends with a married and unloved woman.
Here there are $k$ men and $k$ women and $k$ marriages and $k-1$ affairs.
Hence the BiEGF is $x^ky^{k-1}t^ks^k/k!^2$ ($k \geq 1)$.
But there are $k!^2$ ways of naming (i.e. labelling) 
the persons involved, hence the BiEGF of this case is:
$$
\sum_{k=1}^{\infty} x^ky^{k-1}t^ks^k= {{xts} \over {1-xyts}} \quad .
$$
 
Case IIa: Boris loves Natasha, Natasha is married to Sasha,
Sasha loves Lena, Lena is married to Misha, Misha loves Sofia, ...
 
This case is like the previous case with love ($y$) and marriage ($x$)
interchanged. Hence the BiEGF of this case is:
$$
\sum_{k=1}^{\infty} x^{k-1}y^kt^ks^k= {{yts} \over {1-xyts}} \quad .
$$
 
Case III: Boris is married to Natasha, Natasha loves Andrei,
Andrei is married to Lena, Lena loves Tolya, ...
a chain that starts with a married and unloved
man and ends with a loved and unmarried man.
Here there are $k+1$ men and $k$ women,
and $k$ marriages and $k$ affairs ($k \geq 1$).
Hence the BiEGF is $x^ky^kt^{k+1}s^k/((k+1)!k!)$ ($k \geq 1)$.
But there are $(k+1)!k!$ ways of naming the persons involved, hence
the BiEGF of this case is:
$$
\sum_{k=1}^{\infty} x^ky^kt^{k+1}s^k= {{xyt^2s} \over {1-xyts}} \quad .
$$
 
Case IIIa: Natasha is married to Boris, Boris loves Anna,
Anna is married to Senia, Senia loves Katya, ...
a chain that starts with a married and unloved
woman and ends with a loved and unmarried woman.
 
This is the same as
case III with $t$ and $s$ interchanged. Hence the BiEGF is:
$$
\sum_{k=1}^{\infty} x^ky^kt^{k}s^{k+1}= {{xyts^2} \over {1-xyts}} \quad .
$$
 
Case IV: Boris loves Natasha, Natasha is married to Sasha,
Sasha loves Lena, Lena is married to Tolya, ... , ... 
Marina is married to Boris. Here we have a {\it cycle} of
$k$ married and loved men and $k$ 
married and loved women, $k$ marriages and $k$ affairs.
The BiEGF for each single such scenario is $x^ky^kt^ks^k/k!^2$.
The number of ways of assigning labels is $k!^2/k$
(we have to divide by $k$ because it is a cycle). Hence
the BiEGF of this case is:
$$
\sum_{k=1}^{\infty} {{x^ky^kt^{k}s^{k}} \over {k}}= 
-log(1-xyts) \quad .
$$
 
It follows that the BiEGF for Marriage-Affair profiles,
$$
\sum_{m=0}^{\infty}\sum_{n=0}^{\infty}
{{H_{m,n}(x) H_{m,n}(y) t^ms^n} \over {m!n!}} \quad ,
$$
is equal to
$$
\exp \left ( t+s+{{xts+yts+xyt^2s+xyts^2} \over {1-xyts}} 
-log(1-xyts) \right ) =
$$
$$
(1-xyts)^{-1} \exp \left( {{t+s+xts+yts} \over {1-xyts}} \right )
\quad . \,\, \halmos
$$
 
{\bf Future Directions}
 
It would be worthwhile to generalize the present result to 
find heterosexual analogs of the
multi-variate formulas of Kibble-Slepian and Louck treated
by Foata in [F3]. There, there are several kinds of relationships,
but a cap of two relationships per person. It might be also
interesting to investigate the $r$-sex case.
 
{\bf References}
 
[BG] E. A. Bender and J. R. Goldman, {\it Enumerative uses of
generating functions}, Indiana Univ. Math. J. {\bf 20} (1971),
753-764.
 
[F1] D. Foata, {\it A Combinatorial proof of the Mehler formula},
J. Comb. Th. Ser. A {\bf 24} (1978), 367-378.
 
[F2] D. Foata, {\it ``La S\'erie G\'en\'eratrice Exponenetielle dans 
les Probl\'emes d'\'Enumeration''}, University of Montreal Press,
1971.
 
[F3] D. Foata, {\it Some Hermite polynomial identities and their
combinatorics}, Adv. in Appl. Math. {\bf 2} (1980), 250-259.
 
[FS] D. Foata and M. -P. Sch\"utzenberger, {\it ``Th\'eorie G\'eometrique
des Polyn\^omes Euleriens''}, Lecture Notes in Mathematics {\bf 138},
Springer-Verlag, Berlin, 1970.
 
[W] H. S. Wilf, {\it ``generatingfunctionology''}, 
Academic Press, Boston, 1990.

\bye